\input amstex
\documentstyle{amsppt}
\magnification=\magstep1 

\pagewidth{6.5truein}
\pageheight{9.375truein}
\vcorrection{-0.375truein}

\NoBlackBoxes

\long\def\ignore#1\endignore{#1}

\input xy \xyoption{matrix} \xyoption{arrow} \xyoption{curve} 
\def\edge{\ar@{-}}
\def\dttdar{\ar@{.>}}
\def\dshdar{\ar@{-->}}

\def\bfC{{\bold C}}
\def\bfE{{\bold E}}
\def\bfF{{\bold F}}
\def\bfM{{\bold M}}
\def\bfN{{\bold N}}
\def\bfP{{\bold P}}
\def\bfT{{\bold T}}
\def\bfX{{\bold X}}
\def\bfY{{\bold Y}}
\def\bfZ{{\bold Z}}
\def\bfd{{\bold d}}

\def\A{{\Cal A}}

\def\G{{\Cal G}}

\def\S{{\Cal S}}
\def\U{{\Cal U}}
\def\V{{\Cal V}}

\def\AA{{\Bbb A}}
\def\NN{{\Bbb N}}
\def\ZZ{{\Bbb Z}}

\def\l{\ell}
\def\opE{\operatorname{Ext}^1}

\def\im{\operatorname{Im}}
\def\id{\operatorname{id}}
\def\op{\operatorname{op}}
\def\chr{\operatorname{char}}
\def\Pic{\operatorname{Pic}}
\def\AMod{A{\operatorname{\text{-}Mod}}}
\def\Amod{A{\operatorname{\text{-}mod}}}
\def\BMod{B{\operatorname{\text{-}Mod}}}
\def\Bmod{B{\operatorname{\text{-}mod}}}
\def\ModB{{\operatorname{Mod\text{-}}}B}
\def\modB{{\operatorname{mod\text{-}}}B}
\def\KMod{K{\operatorname{\text{-}Mod}}}
\def\moddash{\operatorname{mod\text{-}}}
\def\dashmod{\operatorname{\text{-}mod}}
\def\Aproj{A{\operatorname{\text{-}proj}}}
\def\add{\operatorname{add}}
\def\End{\operatorname{End}}
\def\Hom{\operatorname{Hom}}
\def\Ext{\operatorname{Ext}}
\def\ker{\operatorname{Ker}}

\def\Out{\operatorname{Out}}
\def\OA{\operatorname{Out}(A)^0}
\def\OB{\operatorname{Out}(B)^0}
\def\Inn{\operatorname{Inn}}
\def\Aut{\operatorname{Aut}}
\def\Pic{\operatorname{Pic}}
\def\DPic{\operatorname{DPic}}
\def\Stab{\operatorname{Stab}}
\def\stab{\operatorname{stab}}
\def\T{T}
\def\d{\operatorname{d}}
\def\can{\operatorname{\frak{can}}}
\def\res{\operatorname{\frak{res}}}
\def\moda#1{\operatorname{Mod}^A_{#1}}
\def\comp{\operatorname{Comp}^A_{\bfd}}
\def\tcomp{T_{X}(\comp)}
\def\GL{\operatorname{GL}}

\def\Stab{\operatorname{Stab}}
\def\Cminus{C^{-}(\Amod)}
\def\Kminus{K^{-}(\Amod)}
\def\Dminus{D^{-}(\Amod)}
\def\Cb{C^b(\Amod)}
\def\Kb{K^b(\Amod)}
\def\Db{D^b(\Amod)}

\def\Orb{\operatorname{\underline{\frak{Orb}}}}
\def\orb{\operatorname{\underline{\frak{orb}}}}

\def\Xtilde{\widetilde{\bfX}}
\def\Ahat{\widehat{A}}
\def\Bhat{\widehat{B}}
\def\Bhat{\widehat{B}}

\def\xtilde{\tilde{\bfX}}

\def\sigbar{\overline{\sigma}}
\def\taubar{\overline{\tau}}

\def\ltensor{\overset {\sssize L} \to\otimes}

\def\AnFu{{\bf 1}}
\def\AuSm{{\bf 2}}
\def\Bor{{\bf 3}}
\def\CPS{{\bf 4}}
\def\DlP{{\bf 5}}
\def\Fro{{\bf 6}}
\def\GuSa{{\bf 7}}
\def\Hap{{\bf 8}}
\def\Kel{{\bf 9}}
\def\Kra{{\bf 10}}
\def\LeMe{{\bf 11}}
\def\Lin{{\bf 12}}
\def\Pol{{\bf 13}}
\def\RicI{{\bf 14}}
\def\Ric{{\bf 15}}
\def\Rou{{\bf 16}}
\def\RoZi{{\bf 17}}
\def\Ver{{\bf 18}}
\def\Voi{{\bf 19}}
\def\Yek{{\bf 20}}

\topmatter

\title Geometry of chain complexes and outer automorphisms under derived
equivalence   \endtitle

\rightheadtext{Geometry of chain complexes and outer automorphisms}

\author Birge Huisgen-Zimmermann and Manuel Saor\'\i n \endauthor

\dedicatory The authors wish to dedicate this paper to Idun Reiten on the
occasion of her sixtieth birthday. \enddedicatory 

\address Department of Mathematics, University of California, Santa Barbara, CA
93106, USA\endaddress

\email birge\@math.ucsb.edu\endemail

\address Departamento de M\'atematicas, Universidad de Murcia, 30100
Espinardo-MU, Spain \endaddress

\email msaorinc\@fcu.um.es \endemail


\abstract The two main theorems proved here are as follows: If $A$ is a finite
dimensional algebra over an algebraically closed field, the identity component of
the algebraic group of outer automorphisms of $A$ is
invariant under derived equivalence. This invariance is obtained as a consequence
of the following generalization of a result of Voigt. Namely, given an appropriate
geometrization
$\comp$ of the family of finite $A$-module complexes with fixed sequence $\bfd$ of
dimensions
  and an ``almost projective'' complex $X\in \comp$, there exists a canonical
vector space embedding 
$$T_{X}(\comp) / T_{X}(G.X) \ \longrightarrow\ \ \Hom_{D^b (\AMod)}(X,
X[1]),$$
where $G$ is the pertinent product of general linear groups acting on $\comp$,
tangent spaces at $X$ are denoted by  $T_X(-)$, and $X$ is identified
with its image in the derived category $D^b (\AMod)$.
\endabstract

\thanks While carrying out this project, the first-named author was partially
supported by an NSF grant, and the second-named author by grants from the DGES of
Spain and the Fundaci\'on `S\'eneca' of Murcia.
The contents of this article were presented at the Conference on Representations
of Algebras at Sa\~o Paolo in July 1999 and at the Workshop on Interactions
between Algebraic Geometry and Noncommutative Algebra at the Mathematical Sciences
Research Institute (Berkeley) in February 2000, by the second- and first-named
authors, respectively.
\endthanks

\endtopmatter

\document

\head{1. Introduction}\endhead

One of our primary goals is to prove that the identity component $\OA$ of the
outer automorphism group of a finite dimensional algebra $A$ is invariant under
derived equivalence; here we assume that the base field $K$ of $A$ is
algebraically closed.  This generalizes a result of Brauer (see
\cite{\Pol}), guaranteeing that $\OA$ is Morita invariant, as well as the next
step beyond Brauer's result which established tilting-cotilting invariance of
$\OA$ (the latter is due to Guil Asensio and the second-named author
\cite{\GuSa}).  On the other hand, note that the full outer automorphism group
$\Out(A) = \Aut(A)/\Inn(A)$ is not even Morita invariant, in general, and a
fortiori, invariance fails for the group $\Aut(A)$ of all algebra automorphisms
of $A$.  Our invariance theorem was independently and with different methods
proved by Rouquier \cite{\Rou}.

In his seminal paper \cite{\RicI}, Rickard extended results of Happel \cite{\Hap}
and Cline-Parshall-Scott \cite{\CPS}, to give an explicit description of the
equivalences between derived categories of module categories (strictly speaking,
between categories derived from the triangulated homotopy categories of chain
complexes of modules).  This characterization opened up the possibilty of
exhibiting invariants under derived equivalence:  Rickard himself showed that the
center and the Grothendieck group of $A$ are among them \cite{\RicI}, as are the
cyclic and Hochschild cohomologies of $A$ (see
\cite{\Kel} and
\cite{\Ric}).   The group $\OA$, which our final theorem below places on this
list, can be viewed as a carrier of homological information as well;  indeed, as
was first observed by Fr\"ohlich in \cite{\Fro},
$\Out(A)$ naturally embeds into the Picard group of $A$, that is, into the group
of (isomorphism types of) Morita self-equivalences of the category of left
$A$-modules.  The invariance status of the full Picard group relative to derived
equivalence is negative, but self-injective algebras which are derived equivalent
share at least the stable Picard group \cite{\Lin}.  To round off the picture, we
point out that the classical Picard group is an object of natural significance in
the context of derived categories, in that it in turn embeds into the `derived
Picard group' $\DPic(A)$, which consists of the isomorphism types of derived
self-equivalences of $A$ induced by tilting complexes in $D^b((A\otimes_K
A^{\op})\text{-Mod})$ (see
\cite{\RoZi}, \cite{\Yek}, and \cite{\LeMe}).  The second of our main results
entails that $\DPic(A)$ contains only finitely many $A$-$A$ tilting
complexes of fixed total dimension which are pairwise non-isomorphic when viewed
as one-sided complexes over $A$. 

This latter result is obtained (in Section 2) as an ingredient of our invariance
proof for
$\OA$.  We describe it in some detail, since it holds substantial interest in its
own right.  Suppose $X$ is a point in the classical variety
$\moda{d}$ of all
$d$-dimensional left
$A$-modules.  In \cite{\Voi}, Voigt exhibited a natural vectorspace monomorphism 
$$T_X(\moda{d})/T_X(\GL_d.X)\ \   \longrightarrow\ \ \Ext^1_{A}(X,X),$$  where
$T_X(\moda{d})$ and $T_X(\GL_d.X)$ denote the tangent spaces at $X$ to $\moda{d}$
and the $\GL_d$-orbit of $X$ in $\moda{d}$, respectively.  Here is a sketch of
our generalization (Theorem 7):  Roughly speaking, it relates tangent spaces of
varieties of bounded finite dimensional complexes over
$A$ to Hom-groups in the derived category
$D^b(\AMod)$, whenever $A$ is finite dimensional over an algebraically closed
field.  More precisely, the geometrization of the
$d$-dimensional modules in the framework of
$\moda{d}$ can be carried over, in the same spirit, to the complexes of the form
$0 \rightarrow X_m
\rightarrow X_{m-1}
\rightarrow
\cdots
\rightarrow X_0 \rightarrow 0$ for fixed $m$, where $X_0, \dots, X_m \in \AMod$
have prescribed dimensions $d_0, \dots, d_m$, respectively.  In the predictable
manner, this leads to a Zariski-closed subset $\comp$ of 
$$\prod_{0 \le i \le m} \Hom_K(A, \End_K(K^{d_i})) \times \prod_{1 \le i
\le m} \Hom_K(K^{d_i}, K^{d_{i-1}}),$$ where $\bold{d} = (d_m, \dots, d_0)$. 
Moreover, the orbits of the (only plausible) conjugation action of
$G = \GL_{d_0}
\times \cdots \times
\GL_{d_m}$ on $\comp$ are in one-to-one correspondence with the isomorphism types
of complexes of the described ilk.  In this setting, the following is true for any
complex $\bfX$ of finite dimensional left
$A$-modules of the above format, with the additional property that all $X_i$,
except for possibly $X_m$, are projective over $A$:  If
$\bfX$ is represented by a point
$X \in \comp$ say, there exists a canonical vector space embedding
$$T_{X}(\comp) / T_{X}(G.X) \ \longrightarrow\ \ \Hom_{D^b (\AMod)}(\bfX,
\bfX[1]).$$    In case all $X_i$ are projective, this embedding is actually an
isomorphism.  Observe that, in
$D^b(\AMod)$, every finite dimensional left
$A$-module coincides with a complex of the addressed ilk  --  to wit, with the
stalk complex of the module concentrated in degree $m$  --  and hence Voigt's
result is retrieved as a special case.  Among the consequences, the one
instrumental in establishing our result concerning
$\OA$ is this:  For any finite sequence of non-negative integers, there are only
finitely many complexes of projectives with this sequence of dimensions such that 
$\Hom_{D^b(\AMod)}(\bold{X}, \bold{X}[1]) =0$.  So, in particular, we see that
$A$ has only finitely many one-sided tilting complexes of prescribed
total dimension.
\smallskip 

Throughout, $A$ will be a finite dimensional algebra over a field $K$. In the
main results of Section 2 and throughout Section 3, $K$ will be algebraically
closed. Moreover,
$\AMod$,
$\text{Mod-}A$,
$\Amod$, and $\moddash A$ will denote the categories of all left/right
$A$-modules, and their full subcategories of finite dimensional left/right
$A$-modules, respectively.

\head{2. The geometry of chain complexes}\endhead

 Suppose that $A$ is generated by 
$a_1=1,a_2, \dots, a_s$ as a $K$-algebra, and recall that the objects in
$\Amod$ of vector space dimension
$d$ are parametrized by the points of the following subvariety of
$s \cdot d^2$-dimensional affine space:
$$\moda{d}= \{ (\A_1, \dots, \A_s) \mid \text{the\ }\A_j\in \End_K(K^d) \text{\
satisfy the same relations as\ } a_1,\dots,a_s \}.$$ This variety comes with a
canonical morphic $\GL_d$-action by conjugation, the orbits of which reflect
basis change. So, if one assigns to each point in $\moda{d}$ the corresponding
left $A$-module, it is clear that the isomorphism types of $d$-dimensional
modules are in one-to-one correspondence with the $\GL_d$-orbits in $\moda{d}$. 
Note that we can identify the points of $\moda{d}$ with points in
$\Hom_K(A,\End_K(K^d))$ by assigning to each sequence $(\A_1, \dots, \A_s) \in
\moda{d}$ the $K$-algebra homomorphism $A
\rightarrow \End_K(K^d)$ which sends $a_j$ to $\A_j$.

We carry this idea over  -- in the same spirit  --  to finite chain complexes
$$ \bfX: \qquad\quad 0\rightarrow X_m @>{\partial_m}>>
X_{m-1}@>{\partial_{m-1}}>> \cdots @>{\partial_1}>> X_0\rightarrow 0$$ in $\Amod$
with fixed sequence 
$$\bfd= (\dim_K X_m, \dots, \dim_K X_0)= (d_m,\dots,d_0)$$ of dimensions.
Accordingly, we define the variety $\comp$ to be the following Zariski
closed-subset of the affine space 
$$\AA^N \ = \  \prod_{0 \le i \le m} \Hom_K(A, \End_K(K^{d_i})) \times \prod_{1
\le i
\le m} \Hom_K(K^{d_i}, K^{d_{i-1}}).$$ Namely,
$$\multline \comp := \{ (\A^m, \dots, \A^0, \partial_m, \dots, \partial_1) \in
\AA^N \mid\ \A^i = (\A_{ij})_{j \le s} \in \moda{d_i},\\
\partial_i \A_{ij}= \A_{i-1,j} \partial_i, \text{\ and\ }
 \partial_{i-1}\partial_i =0 \text{\ for all\ } i \le m \text{\ and\ } j\le s\}.
\endmultline$$ We will also label the points of $\comp$ in the form $(\A_{ij},
\partial_m,
\dots, \partial_1)$, as convenience dictates.  Clearly, the group $\GL_{\bfd}=
\GL_{d_m}
\times\cdots\times \GL_{d_0}$ acts morphically on $\comp$ via
$$(g_m,\dots, g_0)(\A_{ij}, \partial_m, \dots, \partial_1) = (\A_{ij}^{g_i},
g_{m-1}\partial_mg_m^{-1}, \dots, g_0\partial_1g_1^{-1}),$$ where $\A_{ij}^{g_i}$
stands for the conjugate $g_i\A_{ij}g_i^{-1}$.   Moreover, it is straightforward
to check the following 

\proclaim{Observation 1} The natural map from $\comp$ to the collection of
complexes $\bfX$ as above, with fixed sequence $\bfd$ of dimensions, induces a 
one-to-one correspondence between the orbits of the $\GL_{\bfd}$-action on
$\comp$ on one hand, and the isomorphism classes of complexes of the described
format on the other. \qed\endproclaim

An alternate way of geometrizing complexes is, of course, that of viewing
complexes as graded modules over an enlarged algebra; but this approach would be
ill-suited to our purpose of relating the geometry of $\comp$ to the derived
category of $A$. 

In the following, we will denote points of
$\comp$ by capital letters and the corresponding complexes by bold versions of
these letters.  Furthermore, we will, from now on, assume that $a_1=1, \dots,
a_s$ is a
$K$-basis of
$A$.  This provides us with structure constants $c_{jkl} \in K$, arising in the
equalities
$a_ja_k= \sum_{l=1}^s c_{jkl}a_l$.  In order to describe the coordinate ring of
$\comp$, we let
$X_{ij}$ stand for a
$d_i\times d_i$ matrix of variables representing the entries of the matrix
$\A_{ij}$, respectively, and
$Y_i$ for a
$d_{i-1}\times d_i$ matrix of variables representing the entries of the map
$\partial_i$.  When we view a sequence $(X_{ij}, Y_m,
\dots, Y_1)$ as an element of the coordinate ring of $\comp$, our setup shows
these coordinate matrices to be subject to the following equalities; in fact,
$\comp$ is determined by these requirements.

\quad ($\alpha$) $X_{ij}X_{ik}= \sum_{l=1}^s c_{jkl}X_{il}$;

\quad ($\beta$) $Y_iX_{ij}= X_{i-1,j}Y_i$;

\quad ($\gamma$) $Y_{i-1}Y_i= 0$.

In the upcoming lemma, we will derive a convenient explicit description of the
Zariski tangent space $\tcomp$ of $\comp$ at a point
$X$. Recall that a derivation of a left
$A$-module $M$ is a $K$-linear map $\delta: A\rightarrow \End_K(M)$ such that
$\delta(ab)x= a(\delta(b)x)+ \delta(a)(bx)$ for all $a,b\in A$ and $x\in M$.

\proclaim{Lemma 2} Given a point $X= (\A^m, \dots, \A^0, \partial_m, \dots,
\partial_1) \in \comp$, the tangent space $\tcomp$ of $\comp$ at $X$ consists
precisely of those sequences
$(\delta_m, \dots, \delta_0, \sigma_m, \dots, \sigma_1)$ in 
$$\AA^N \ = \  \prod_{0 \le i
\le m} \Hom_K(A, \End_K(K^{d_i})) \times \prod_{1 \le i
\le m} \Hom_K(K^{d_i}, K^{d_{i-1}})$$ which satisfy the following conditions:

{\rm (a)} Each $\delta_i : A\rightarrow \End_K(K^{d_i})$ is a derivation of the
left $A$-module $X_i$ determined by $\A^i = (\A_{ij})_{j\le s}$.

{\rm (b)} $\sigma_i \A_{ij}+ \partial_i \delta_i(a_j)= \A_{i-1,j}
\sigma_i+ \delta_{i-1}(a_j) \partial_i$  for all $0 \le i \le m$ and $1 \le j \le
s$; here we identify $\A_{ij}$ with the left multiplication of $X_i = K^{d_i}$ by
the element
$a_j$ of $A$.

{\rm (c)} $\sigma_{i-1} \partial_i+ \partial_{i-1} \sigma_i =0$ for all
$1 \le i \le m$. \endproclaim

\demo{Proof}  Note that each of the $\delta_i$ is determined by the values
$(\delta_i(a_j))_{j\le s}$.  Hence the defining conditions
$(\alpha)-(\gamma)$ of $\comp$ yield the following equations pinning down
$\tcomp$; they are given in the variables $w_{ij}$ representing the maps
$\delta_i(a_j)$, and $v_i$ representing the $K$-endomorphisms
$\sigma_i: K^{d_i} \rightarrow K^{d_{i-1}}$, respectively:

\quad (a) $w_{ij} \A_{ik} + \A_{ij} w_{ik} = \sum_{l=0}^s c_{jkl}w_{il}$,

\quad (b) $v_i \A_{ij} + \partial_i w_{ij} = w_{i-1,j} \partial_i +
\A_{i-1,j} v_i$,

\quad (c) $v_{i-1} \partial_i + \partial_{i-1} v_i= 0$,

\noindent for all eligible indices $i,j,k$.  But these conditions are clearly
tantamount to the ones listed in our claim. \qed 
\enddemo

Next we follow the traditional road of mapping the tangent space of a module
variety at a point $x$ to the group of self-extensions of the module represented
by x.  Namely, we fix
$X \in \comp$, and assign to any point $(\delta_m, \dots, \delta_0, \sigma_m,
\dots,
\sigma_1)$ in the tangent space $\tcomp$ a short exact sequence $0
\rightarrow \bfX
\rightarrow \bfZ \rightarrow \bfX \rightarrow 0$ of complexes as follows. 

\definition{Remark 3} Let $X \in \comp$ as before.

{\rm (1)} Condition {\rm (a)} of Lemma 2 implies that the following is a left
$A$-module structure on $Z_i = K^{d_i}\oplus K^{d_i}$:  Namely, we define left
multiplication of $Z_i$ by $a \in A$ as left multiplication by $\left(
\smallmatrix
\A^i(a) &\delta_i(a)\\ 0 &\A^i(a) \endsmallmatrix \right)$, viewing the elements
$Z_i$ as column vectors.

{\rm (2)} Condition (c) of Lemma 2 tells us that we can supplement the $Z_i$ to a
chain complex
$\bfZ$ of $K$-vector spaces by introducing the differentials
$\partial_i^{\bfZ}=\left(
\smallmatrix
\partial _i & \sigma _i\\ 0 & \partial _i \endsmallmatrix \right)$.

{\rm (3)}  Condition (b) of Lemma 2, finally, shows that the complex $\bfZ$ just
defined is, in fact, a complex of $A$-modules. 

{\rm (4)}  Let $0\rightarrow X_i\rightarrow Z_i\rightarrow X_i\rightarrow 0$ be
the canonical sequences, given by injection of $X_i$ into the first component of
$Z_i$, followed by projection onto the second component.  These short exact
sequences in $\Amod$ are compatible with the differentials, thus yielding a short
exact sequence
$0\rightarrow\bfX\rightarrow\bfZ\rightarrow\bfX\rightarrow 0$ of chain complexes.
$\qed$ 
\enddefinition

Now suppose that $\bfX , \bfY$ are chain complexes which are bounded above, but
not necessarily finite.  We will consider the group $\opE(\bfX,\bfY)$ of
(equivalence classes) of extensions $0 \rightarrow \bfY \rightarrow \bfE
\rightarrow \bfX \rightarrow 0$ in the category of complexes of $A$-modules.  In
analogy with  Ext-groups of $A$-modules, 
$\opE(\bfX,\bfY)$ actually carries a $K$-vector space structure.  

For any point $X \in \comp$, Remark 3 provides us with a
$K$-linear map
$$\chi: \tcomp \rightarrow \opE(\bfX,\bfX),$$
 sending an element of $\tcomp$ to the class of the self-extension $0 \rightarrow
\bfX \rightarrow \bfZ \rightarrow \bfX \rightarrow 0$ of $\bfX$ described above. 
The proof of the following proposition is quite similar to that of Voigt's Lemma
(\cite{\Voi, Chap\.2, Section 3.4}; see also
\cite{\DlP, Theorem 1.6}).  Since Voigt's original argument (couched in German)
carries over smoothly to our setting, we include only a brief outline (for the
sake of the reader not familiar with German).  Recall that an extension $0
\rightarrow \bfY \rightarrow \bfE \rightarrow
\bfX \rightarrow 0$ is {\it semisplit} in the sense of Verdier \cite{\Ver, p\.
272-273} if, in each degree $i$, the pertinent short exact sequence of
$A$-modules splits.  

\proclaim{Proposition 4}  Suppose that $K$ is an algebraically closed field. 
Then, given any $X \in \comp$, the kernel of the map $\chi$ equals the tangent
space of
$G.X$ at
$X$.  In other words,
$\chi$ induces a vectorspace embedding
$$\tcomp / T_{X}(G.X)  \rightarrow \opE(\bfX,\bfX).$$ Moreover, all semisplit
self-extensions of $\bfX$ belong to the image of $\chi$.  
\endproclaim

\demo{Proof} The canonical morphism $G \rightarrow G.X$ is separable (since the
stabilizer subgroup of $X$ in $G$ arises as the solution set of a system of linear
equations), i.e., this map has a surjective differential   
$$T_e(G) \rightarrow T_{X}(G.X) \ \subseteq \ \tcomp.$$  In the description of the
tangent space
$\tcomp$ given in Lemma 2, this identifies $T_X(G.X)$ as being the following
subset of
$\tcomp$: Namely, the set of those elements $(\delta_m, \dots, \delta_0,
\sigma_m, \dots, \sigma_1)$ in $\tcomp$ for which there is a sequence $t = (t_m,
\dots, t_0) \in \prod_{i=0}^m \End_K(K^{d_i})$ such that the following equalities
hold for all eligible $i$ and $j$:

$\bullet$  $\delta_i(a_j) = t_i \A_{ij} - \A_{ij}t_i$, and

$\bullet$  $\sigma_i = t_{i-1} \partial_i - \partial_i t_i$.

\noindent (These equalities can be gleaned directly from Voigt's argument for
modules by identifying the complex $\bfX$ with a graded module $\bigoplus_{i=0}^m
X_i$ over the finite dimensional algebra $A[\partial_1, \dots,
\partial_m]/(\partial_1^2, \dots,
\partial_m^2)$;  here the operation of $A$ on the direct sum is extended via
$\partial_i.(x_m, \dots, x_0) = \partial_i(x_i) \in X_{i-1}$ for $1 \le i \le
m$.)  

It is readily verified that the described subset of $\tcomp$ coincides with the
kernel of
$\chi$. Indeed, splitness of the extension $0 \rightarrow \bfX \rightarrow \bfZ
\rightarrow \bfX
\rightarrow 0$ of complexes coming with the point $(\delta_m, \dots, \delta_0,
\sigma_m,
\dots, \sigma_1)$ of $\tcomp$ is equivalent to the existence of a section for the
epimorphism
$\bfZ
\rightarrow \bfX$.  Since the latter map just projects the
$Z_i = X_i \oplus X_i$ onto the second components, such sections are precisely
the chain maps $\bfX \rightarrow \bfZ$ of the form $t' = (t_i,1_{X_i})_{0 \le i
\le m}$, where the family
$(t_i)$ belongs to $\prod_{i=0}^m \End_A(X_i)$.  But that $t_i': X_i
\rightarrow Z_i$ be an 
$A$-module homomorphism is clearly tantamount to the requirement that $t_i
\A_{ij} =
\A_{ij}t_i +
\delta_i(a_j)$, while the stipulation that a family $(t_i) \in \prod_{i=0}^m
\Hom_K(K^{d_i})$ be compatible with the differentials of $\bfX$ and $\bfZ$
amounts to the equality $t_{i-1} \partial_i =
\partial_i t_i  +  \sigma_i$.  This completes the proof of the fact that
$\T_X(G.X)$ coincides with the kernel of $\chi$.

For the final comment, note that the differential of any semisplit self-extension
$0 \rightarrow \bfX
\rightarrow \bfZ \rightarrow \bfX \rightarrow 0$ of $\bfX$, with splitting $Z_i =
X_i \oplus X_i$ in degree
$i$, has the form $\partial_i^{\bfZ}=\left(
\smallmatrix
\partial _i & \sigma _i\\ 0 & \partial _i \endsmallmatrix \right)$ relative to
such a splitting, where $\partial_i$ is the $i$-th differential of $\bfX$ and
$\sigma_i \in \Hom_A(X_i, X_{i-1})$ satisfies condition (c) of Lemma 2. 
Consequently, Remark 3 shows that the semisplit self-extensions of $\bfX$ are
precisely the images under
$\chi$ of the vectors $(0, \dots, 0, \sigma_m, \dots, \sigma_1)$ in the tangent
space $\tcomp$.   \qed 
\enddemo

Following Verdier's conventions in \cite{\Ver}, we denote the category of right
bounded complexes of finite dimensional left $A$-modules by $\Cminus$, its
quotient category modulo homotopy by $\Kminus$, and the corresponding derived
category by $\Dminus$.  Analogously, $\Cb$, $\Kb$, and $\Db$ stand for the
category of bounded complexes, for the quotient category modulo homotopy, and the
pertinent derived category, respectively.  
  According to
\cite{\Ver, pp\. 294-295}, every element
$0\rightarrow\bfY\rightarrow\bfE\rightarrow\bfX\rightarrow 0$ of
$\opE(\bfX,\bfY)$ gives rise to a distinguished triangle
$\bfY\rightarrow\bfE \rightarrow\bfX\rightarrow\bfY[1]$  in $\Dminus$.  Here
$\bfY[1]$ denotes the shifted complex which carries $Y_{n-1}$ in the slot labeled
$n$. This provides us with a
$K$-linear map
$\xi :\opE(\bfX,\bfY)\longrightarrow
\Hom_{\Dminus}(\bfX ,\bfY[1])$ which assigns to an extension as above the
connecting morphism
$\bfX\rightarrow\bfY[1]$ of the corresponding triangle. Our next intermediate
goal is to scrutinize this map for complexes
$\bfX$ which are `close' to being projective.

\definition{Definition} A complex $\bfX : \ \cdots \rightarrow X_{i+1}
\rightarrow X_i \rightarrow
\cdots$ in $\Cminus$ is called {\it projective} in case all of the terms $X_i$ are
projective $A$-modules.   

Now suppose that  $\bfX \in \Cb$.  If $\bfX$ is nonzero, then the largest integer
$m$ with $X_m \ne 0$ will be referred to as the {\it left degree\/} of $\bfX$. We
will call $\bfX$ {\it almost projective\/} in case either
$\bfX$ is projective or else $X_m$ is the only non-projective term of the
complex.   
\enddefinition

It is clear that, under the natural map $\Cb \rightarrow \Db$,  every object in
the derived category is represented by an almost projective complex.  

One of the assets of almost projective complexes is a convenient lifting property
for maps in the derived category:  If $\bfX$ and $\bfY$ are both almost
projective with coinciding left degree, then any map in $\Hom_{\Db}(\bfX, \bfY)$
lifts to a chain map in
$\Hom_{\Cb}(\bfX, \bfY)$:  Indeed, we may restrict our attention to bounded
almost projective complexes for that purpose.  So suppose that 
$\bfX$ is of the form
$$0 \rightarrow X_m @>\partial_m>> X_{m-1} @>\partial_{m-1}>> \dots
@>\partial_1>> X_0
\rightarrow 0,$$ where all of the left $A$-modules $X_i$, except for $X_m$
possibly, are projective; write $\bfY$ in the same format, with terms $Y_i$.  If
we let $\cdots \rightarrow P_{m+1} @>p_{m+1}>> P_m @>p_m>> X_m
\rightarrow 0$ be a projective resolution of $X_m$, the image of $\bfX$ in
$\Dminus$ clearly coincides with that of the projective complex $\bfX'$, given by
$$\dots P_{m+1} @>p_{m+1}>> P_m @>\partial_m p_m>> X_{m-1} @>\partial_{m-1}>>
\dots @>\partial_1>> X_0 \rightarrow 0.$$   Analogously, define a projective
complex $\bfY'$ based on $\bfY$.  Now any homomorhism in $\Hom_{\Dminus}(\bfX',
\bfY')$ is induced by a chain map
$\eta \in \Hom_{\Cminus}(\bfX', \bfY')$.  But $p_m$ is the cokernel of $p_{m+1}$,
and therefore $\eta$ induces a chain map in
$\Hom_{\Cb}(\bfX, \bfY)$ having the same image as $\eta$ in $\Db$.  

Again suppose that $\bfX \in \Cb$ is almost projective.  Along the line suggested
by the previous paragraph, we associate another almost projective complex
$\bfX^{(1)}$ to $\bfX$, by potentially extending $\bfX$ to the left as follows: 
Set $\bfX^{(1)} = \bfX$ if
$\bfX$ is projective.  Otherwise, we denote the left degree of $\bfX$ by $m$,
let  $p: X_m^{(1)} \rightarrow X_m$ be a projective cover of $X_m$ with kernel
$\partial^{\bfX^{(1)}}_{m+1}: X_{m+1}^{(1)}
\rightarrow X_m^{(1)}$, and define
$\bfX^{(1)}$ to be the complex
$$0 \rightarrow X_{m+1}^{(1)} \rightarrow X_m^{(1)} \rightarrow X_{m-1}
\rightarrow \cdots \rightarrow X_0 \rightarrow 0,$$  where the $m$-th
differential of the complex $\bfX^{(1)}$ is $\partial^{\bfX^{(1)}}_m =
\partial^{\bfX}_m\, p$.  Clearly,
$\bfX^{(1)}$ comes with a natural map $f: \bfX^{(1)} \rightarrow \bfX$ of chain
complexes, given by $f_m = p$, $f_i = 0$ for $i \ge m+1$, and $f_i = \id_{X_i}$
for $i
\le m-1$.  By setting
$\bfX^{(0)} = \bfX$ and continuing inductively  via $\bfX^{(r+1)} =
\bigl(\bfX^{(r)}\bigr)^{(1)}$, we thus obtain a sequence 
$$\cdots \rightarrow \bfX^{(r)} \rightarrow \bfX^{(r-1)} \rightarrow
\cdots \rightarrow \bfX^{(1)} \rightarrow \bfX^{(0)}$$ of chain maps in
$\Cminus$, which, in turn, induces a sequence       
$$\opE(\bfX,\bfY)\rightarrow\opE(\bfX^{(1)},\bfY)\rightarrow\opE(\bfX^{(2)},
\bfY) \rightarrow \cdots $$ of $K$-linear maps for each $\bfY \in \Cminus$.  Note
that Verdier's map
$\xi :\opE(\bfX,\bfY)\rightarrow
\Hom_{\Dminus}(\bfX,\bfY[1])$ factors through these vectorspace homomorphisms,
since all of the canonical maps
$\bfX^{(r)} \rightarrow \bfX^{(r-1)}$ become isomorphisms in the homotopy category
$\Kminus$.

We begin with an auxiliary point which sometimes allows us to replace $\bfX$ by
$\bfX^{(1)}$ in the first argument of $\opE$ without penalty. 

\proclaim{Lemma 5} Let $\bfX, \bfY \in \Cb$ be nonzero almost projective
complexes whose left degrees coincide. Then the canonical K-linear map $\opE(\bfX,
\bfY)\rightarrow\opE(\bfX^{(1)}, \bfY)$ is bijective.
\endproclaim

\demo{Proof}  The coinciding left degree of $\bfX$ and $\bfY$ is denoted by $m$.
Moreover, we write
$\bfZ$ for $\bfX^{(1)}$, and $f$ for the canonical chain map $\bfZ \rightarrow 
\bfX$; this means, in particular, that $f_m: Z_m \rightarrow X_m$ is a projective
cover, and $f_i = \id_{X_i}$ for $i \le m-1$.  The $K$-linear map addressed in
our claim, finally, is denoted by $h: \opE(\bfX,\bfY)
\rightarrow \opE(\bfZ,\bfY)$.

We start by proving surjectivity of $h$.  Suppose
$$\phi: 0 \rightarrow \bfY @>{\rho}>> \bfF @>{\sigma}>> \bfZ \rightarrow 0$$
represents a class in $\opE(\bfZ,\bfY)$.  To construct a preimage
$$\epsilon: 0 \rightarrow  \bfY @>{\mu}>> \bfE @>{\nu}>> \bfX
\rightarrow 0$$ of (the class of) $\phi$ under $h$, we start by setting $E_i = 0$
for $i \ge m+1$ and $E_i = F_i$ for $i \le m-1$; the pertinent components of the
differential $\partial^{\bfE}$ and those of the chain maps $\mu, \nu$ are as
follows:
$\partial^{\bfE}_i = \partial^{\bfF}_i$, $\mu_i = \rho_i$, and $\nu_i =
\sigma_i$ for $i \le m-1$, the definitions to the left of the $m$-th position
being obvious.  The
$A$-module
$E_m$ and the corresponding differential of $\bfE$ are defined by the
requirements that
$g_m: F_m
\rightarrow E_m$ be the cokernel of $\partial^{\bfF}_{m+1}$, and
$\partial^{\bfE}_m$ be the unique map in $\Hom_A(E_m,F_{m-1})$ such that
$\partial^{\bfF}_m = \partial^{\bfE}_m\, g_m$.  Moreover, set $\mu_m =
g_m\rho_m$, and let $\nu_m: E_m \rightarrow X_m$ be the unique homomorphism with
$f_m \sigma_m = \nu_mg_m$; the map $\nu_m$ is, once more, furnished by the
universal property of $g_m$.  Clearly, $\mu: \bfX \rightarrow
\bfE$ is a chain map.  The same is true for $\nu: \bfE \rightarrow \bfX$; indeed,
a straightforward diagram chase yields $\partial^{\bfX}_m\, \nu_m =
\nu_{m-1} \partial^{\bfE}_m$.  It is now routine to check that $\epsilon$ yields
a preimage of $\phi$ under $h$.

To verify injectivity of $h$, suppose that the class of $\epsilon: 0
\rightarrow \bfY @>{\mu}>> \bfE @>{\nu}>> \bfX \rightarrow 0$ is mapped to zero
by $h$.  This means that the upper row of the following pullback diagram of
complexes splits:

\ignore
$$\xymatrixrowsep{2pc}\xymatrixcolsep{4pc}
\xymatrix{ 0 \ar[r] &\bfY \ar[r]^{\rho} \ar@{=}[d] &\bfF \ar[r]^{\sigma} \ar[d]^g
&\bfZ 
\ar[r] \ar[d]^f &0\\ 0 \ar[r] &\bfY \ar[r]_{\mu} &\bfE \ar[r]_{\nu} &\bfX \ar[r]
&0 }$$
\endignore

\noindent In order to show that $\epsilon$ is trivial, i.e., that the lower row
splits as well, let $\tau = (\tau_n)_{n \in \NN}$ be a section for $\sigma$.  To
construct a section $\pi$ for $\nu$, we start by setting $\pi_i = 0$ for $i
\ge m+1$, and $\pi_i = g_i \tau_i$ for $i \le m-1$.  This setup is recorded in
the following diagram:  

\ignore
$$\xymatrixrowsep{2.25pc}\xymatrixcolsep{1.75pc}
\xymatrix{
 &F_{m+1} \ar[rr] \ar[dl]_{\sigma_{m+1}} \ar[dd]|\hole_(0.67){g_{m+1}} &&F_m
\ar[rr]
\ar[dl]_{\sigma_m} \ar[dd]|\hole_(0.67){g_m} &&F_{m-1} \ar[rr]
\ar[dl]_{\sigma_{m-1}}
\ar[dd]|\hole_(0.67){g_{m-1}} &&F_{m-2}  \ar[dl]_{\sigma_{m-2}}
\ar[dd]_(0.67){g_{m-2}}\\ Z_{m+1} \ar[dd]_(0.33){f_{m+1}} \ar[rr]
\ar@/_/[ur]_{\tau_{m+1}} &&Z_m
\ar[dd]_(0.33){f_m} \ar[rr] \ar@/_/[ur]_{\tau_m} &&Z_{m-1} \ar@{=}[dd]
\ar[rr] \ar@/_/[ur]_{\tau_{m-1}} &&Z_{m-2} \ar@{=}[dd]
\ar@/_/[ur]_{\tau_{m-2}}\\
 &0 \ar[dl]_{\nu_{m+1}} \ar[rr]|{\hole\hole} &&E_m \ar[dl]_{\nu_m}
\ar[rr]|{\hole\hole} &&E_{m-1} \ar[dl]_{\nu_{m-1}}\ar[rr]|\hole &&E_{m-2}
\ar[dl]_{\nu_{m-2}}\\ 0 \ar[rr] \dttdar@/_/[ur]_{\pi_{m+1}} &&X_m \ar[rr]
&&X_{m-1} \ar[rr]
\dttdar@/_/[ur]_{\pi_{m-1}} &&X_{m-2} \dttdar@/_/[ur]_{\pi_{m-2}} 
 }$$
\endignore

\noindent Note that all of the squares in this diagram commute, so that our task
is reduced to finding a section
$\pi_m$ for $\nu_m$ such that the two bottom squares adjacent to
$\pi_m$ commute.  This is automatic for the left-hand square, irrespective of our
choice of
$\pi_m$.  To construct $\pi_m$ with the property that $\partial_m^{\bfE} \pi_m =
\pi_{m-1}
\partial_m^{\bfX}$, we recall that
$f_m$ is the cokernel of the differential $\partial^{\bfZ}_{m+1}$.  This provides
us with a unique homomorphism $\pi_m: X_m \rightarrow E_m$ having the property
that $g_m
\tau_m = \pi_mf_m$, for $g_m \tau_m \partial^{\bfZ}_{m+1} \sigma_{m+1} =
\partial^{\bfE}_{m+1} g_{m+1} = 0$ yields $g_m \tau_m \partial^{\bfZ}_{m+1} =
0$.  It is now routine to check that
$\pi = (\pi_i)_{i
\in \ZZ}$ is a section for $\nu$, as required.  This completes the argument.
$\qed$
\enddemo 

Modules $M,N \in \Amod$ can of course be viewed as stalk complexes $\bfM, \bfN$
concentrated in degree
$0$. In this situation, Verdier's map provides us with isomorphisms
$$\opE(\bfM,\bfN) \cong \Hom_{\Dminus}(\bfM,\bfN[1]).$$   The following lemma
shows that the good behavior of Verdier's map extends, at least partially, to
more general complexes, which provides the missing link towards the main result
of this section.

\proclaim{Lemma 6}  Let $\bfX$, $\bfY$ be nonzero complexes in $C^-(\Amod)$. 
Then Verdier's map 
$$\xi:\opE(\bfX,\bfY)\longrightarrow \Hom_{D^-(\Amod)}(\bfX,\bfY[1])$$ is
injective if either

{\rm (1)} $\bfX$ is a projective complex, or else

{\rm (2)} $\bfX$ and $\bfY$ are bounded almost projective complexes such that the
left degree of $\bfX$ is larger than or equal to the left degree of
$\bfY$.

Moreover, if $\bfX$ is as under {\rm (1)} or {\rm (2)} and $\bfY$ is projective,
then $\xi$ is bijective.
\endproclaim
 
\demo{Proof}  We show injectivity of $\xi$ simultaneously for (1) and (2).  In
case (2), we denote by $k$ the left degree of $\bfX$, by $m$ that of
$\bfY$, and observe that Lemma 5 permits us to derive the claim for $k = m$ from
that for $k = m+1$.  In other words, we may, w.l.o.g., assume that either both
$\bfX$ and $\bfY$ are projective or both are bounded almost projective with $k \ge
m+1$.   Note that, in this situation, any extension $\epsilon: 0\rightarrow \bfY
@>{\mu}>>
\bfE @>{\nu}>> \bfX\rightarrow 0$ has the property that $0 \rightarrow Y_i
@>{\mu_i}>> E_i @>{\nu_i}>> X_i \rightarrow 0$ splits  (i.e., the extension
$\epsilon$ is semi-scind\'{e}e in the sense of
\cite{\Ver, pp\. 272-273}).  For each $i \in \ZZ$, let $\pi'_i$ be a section of
$\nu_i$.  We will see that the assumption `$\xi(\epsilon) = 0$' allows us to
adjust the family $(\pi_i')$ so as to furnish a section for $\nu$. To that end,
suppose that  
$\xi(\epsilon)$ is represented by a chain map $\phi :\bfX\rightarrow\bfY[1]$. 
Letting
$\bfP\rightarrow\bfX$ be the projective resolution of $\bfX$ in $\Kminus$ and
$\tilde{\phi}:\bfP\rightarrow\bfY[1]$ its composition with $\phi$, the morphism
$\tilde{\phi}$ in $\Kminus$ also represents $\xi(\epsilon)$ in
$\Dminus$, in view of the canonical isomorphism $\Hom_{\Kminus}(\bfX,\bfY[1])
\cong
\Hom_{\Kminus}(\bfP,\bfY[1])$ (see \cite{\Ver, p\. 299}).  We further observe
that 
$\xi(\epsilon)=0$ if and only if $\tilde{\phi}$ is trivial in $\Kminus$, for the
quotient map $\Hom_{\Kminus}(\bfP,\bfY[1])\rightarrow
\Hom_{\Dminus}(\bfP,\bfY[1])$ is injective.  Our hypothesis on the left
degrees of
$\bfX$ and $\bfY$ guarantees that
$\tilde{\phi}$ is nullhomotopic precisely when $\phi$ is nullhomotopic.
Consequently, the assumption that $\xi(\epsilon)$ be zero supplies us with a
homotopy $\psi$ from
$\phi$ to zero.   We deduce that
$\pi = (\pi_n' -\mu_n\psi_n)_{n\in\ZZ}$ is a chain map
$\bfX\rightarrow\bfE$ such that $\nu \pi =1_{\bfX}$, which shows that $\epsilon$
is indeed trivial.

Now suppose that $\bfY$ is projective.  In that case, we have isomorphisms
$$\Hom_{\Dminus}(\bfX, \bfY[1]) \ \cong \
\Hom_{\Kminus}(\bfP,\bfY[1]) \ \cong \ \Hom_{\Kminus}(\bfX, \bfY[1]),$$ where
$\bfP$ is a projective cover of $\bfX$.  Hence all we need to check is that the
homotopy class of any chain map $\psi: \bfX \rightarrow \bfY[1]$ occurs in the
image of the map $\Ext^1(\bfX,\bfY) \rightarrow \Hom_{\Kminus}(\bfX,\bfY[1])$
induced by
$\xi$.  But if $\bfC$ is  the mapping cone of $\psi$, then clearly the class of
the extension $0 \rightarrow \bfY \rightarrow \bfC[-1] \rightarrow \bfX
\rightarrow 0$ in
$\Ext^1(\bfX, \bfY)$ is such a preimage of $\psi$. 
\qed\enddemo

From the preceding lemmas it is now easy to glean

\proclaim{Theorem 7} Suppose that $K$ is algebraically closed, and let $\bfX: 0
\rightarrow X_m \rightarrow \dots \rightarrow X_0 \rightarrow 0$ be an almost
projective complex of finite-dimensional A-modules, represented by a point
$X$ in $\comp$, where $\bfd = (\dim X_m, \dots, \dim X_0)$.  Moreover, let
$\tcomp$ and
$T_{X}(G.X)$ be the tangent spaces at $X$ of $\comp$ and
$G.X$, respectively.  Then there exists a canonical $K$-linear embedding
$$\tcomp /T_{X}(G.X)  \longrightarrow \Hom_{\Db}(\bfX,\bfX[1]).$$ 

In case $\bfX$ is a projective complex, this embedding is an isomorphism.
\endproclaim 

\demo{Proof}  We compose the vectorspace homomorphism 
$\chi :\tcomp \rightarrow\opE(\bfX ,\bfX)$ introduced after Remark 3 with
Verdier's map
$\xi: \opE(\bfX,\bfX)\rightarrow
\Hom_{\Dminus}(\bfX,\bfX[1])$ to obtain a vectorspace homomorphism
$\eta: \tcomp \rightarrow \Hom_{D^-(\Amod)}(\bfX,\bfX[1])$.  Observe that the
kernel of
$\eta$ coincides with $T_{X}(G.X)$, since $\ker(\chi) = T_{X}(G.X)$ by
Proposition 4, while $\xi$ is an injection by Lemma 6.  Consequently, $\eta$
induces an embedding of
$\tcomp / T_{X}(G.X)$ into
$\Hom_{D^-(\Amod)}(\bfX,\bfX[1])$ as required.  

For the final claim, let $\bfX$ be projective.  Then $\xi \chi$ is surjective,
since $\chi$ is surjective by the last part of Proposition 4, and $\xi$ is
surjective by the last statement of Lemma 6.
\qed
\enddemo

As in the case of the variety $\moda{d}$ of $d$-dimensional $A$-modules, Theorem
7 entails a number of interesting consequences. We just point out two of them,
the second of which we will need in the sequel. Since the proofs are analogous to
those for the module-theoretic counterparts (see \cite{\DlP}), we leave them to
the reader.

\proclaim{Corollary 8}  Suppose that $K$ is algebraically closed and $\bfX\in
\Cb$ a bounded almost projective complex such that
$\Hom_{\Db}(\bfX,\bfX[1])= 0$.  If $X \in \comp$ is a point representing
$\bfX$, then the G-orbit of $X$ is open in
$\comp$.
\endproclaim

\proclaim{Corollary 9}  Again suppose that $K$ is algebraically closed. Given any
finite sequence $\bfd = (d_m, \dots d_0)$ of natural numbers, there are  --  up
to isomorphism in $\Cb$  --  only finitely many almost projective complexes
$\bfX$ of A-modules satisfying the following two conditions:

{\rm (1)} $\Hom_{\Db}(\bfX,\bfX[1])=0$, and

{\rm (2)} $\dim_KX_n=d_n$, for all $n\in\ZZ$.
\endproclaim

\head{3. Faithfully balanced two-sided complexes and invariance of $\OA$}
\endhead

Several of the ideas underlying this section have precursors which were developed
in a paper by Guil-Asensio and the second-named author \cite{\GuSa}.  Major
portions of our auxiliary results can be readily generalized to algebras over
commutative rings; in particular, this is true for Proposition 10 below.  To
facilitate the reading with a unified blanket hypothesis, we will, however,
continue to assume that $A$ is a finite dimensional algebra over an algebraically
closed field
$K$.  

Of course, we can view a chain complex
$\bfX
\in C(\AMod)$ as a pair $(\bfX,\lambda)$, where
$\bfX\in C(\KMod)$ and $\lambda : A\rightarrow \End_{C(\KMod)}(\bfX)$ is a
K-algebra homomorphism.  Accordingly, given another finite dimensional
$K$-algebra $B$, an $A$-$B$-bimodule complex $_A\bfX_B$ amounts to a chain
complex $\bfX$ of $K$-spaces, combined with an algebra homomorphism $\tau:
A\otimes_K B^{\op}\rightarrow
\End_{C(\KMod)}(\bfX)$.  We will call a bimodule complex $_A\bfX_B$ almost
projective if both $_A\bfX$ and $\bfX_B$ are {\it almost projective} in the sense
of Section 2; recall that, in particular, this means that $\bfX$ is a right
bounded complex of finite dimensional modules.  By adapting standard terminology
for modules to this context (see, e.g.,
\cite{\AnFu}), we moreover define:

\definition{Definition} A  bimodule complex $_A\bfX_B$ will be called {\it
faithfully balanced} provided that both of the canonical K-algebra homomorphisms 
$$\lambda = \tau(-,1) : A\rightarrow \End_{C(\ModB)}(\bfX_B)
\qquad\text{and}\qquad \rho = \tau(1,-) : B\rightarrow
\End_{C(\AMod)}(_A\bfX)^{\op}$$ are isomorphisms.  
\enddefinition

Whenever we let endomorphisms of a complex $\bfY$ of left $A$-modules act on the
right of
$\bfY$, we record this by referring to them as maps in
$\End_{C(\Amod)}(\bfY)^{\op}$. The existence of a faithfully balanced bimodule
complex
$_A\bfX_B$ entails a fairly tight connection between the algebras $A$ and $B$, as
the following two results suggest.

\proclaim{Proposition 10} Whenever there exists a faithfully balanced bimodule
complex $_A\bfX_B$, the centers of A and B are isomorphic. 
\endproclaim 

\demo{Proof}  Set $E = \End_{C(\KMod)}(\bfX)$. The canonical embedding of
algebras  $\lambda : A\rightarrow E$ and  $\rho : B^{\op}\rightarrow E$ take the
centers of A and B to the same subalgebra $E'$ of $E$, namely to $E' =
\End_{C(\AMod\text{-}B)}(\bfX)$.
\qed\enddemo

As in Section 1, we denote by $\Out(A)$ the group of outer automorphisms of $A$,
and by $\OA$ its identity component.   It is well-known (cf. \cite{\Fro}, Theorem
1) that $\Out(A)$ embeds naturally into the Picard group $\Pic(A)$, which can be
thought of as the group of isomorphism types of Morita self-equivalences of
$\AMod$.  Clearly, any such Morita self-equivalence induces category equivalences
$C(\AMod)\rightarrow C(\AMod)$, $\Cb \rightarrow \Cb$, and
$D(\AMod)\rightarrow D(\AMod)$, $\Db \rightarrow \Db$.  This provides us with a
natural action of $\Pic(A)$, and consequently also of
$\Out(A)$, on sets of isomorphism classes of objects in
$C(\AMod)$ (resp\. $D(\AMod)$). Given any complex $\bfX\in C(\AMod)$, we will
accordingly consider its $\Out$-{\it orbit}
$\{^{\sigma}\bfX\in C(\AMod):$ $\sigbar\in \Out(A)\}$ in
$C(\AMod)$, and analogously that in $D(\AMod)$, as well as its $\Out$-{\it
stabilizer}
$\Out^{\bfX}(A)=\{\sigbar\in \Out(A): {}^{\sigma}\bfX\cong\bfX\ \text{in\ }
D(\Amod)\}$ in the derived category.  As we are dealing with bimodule complexes
$_A\bfX_B$, it will promote orientation to denote the actions of $\Out(A)$ and
$\Out(B)$ by means of left and right superscripts, respectively.     

\proclaim{Proposition 11} Suppose that $_A\bfX_B$ is a bounded, faithfully
balanced bimodule complex of finite dimensional modules.  If the
$\Out$-orbits of
$_A\bfX$ and $\bfX_B$ in $\Cb$ and $C^b(\Bmod)$  are finite, then the groups
$\OA$ and
$\OB$ are isomorphic. In particular, this is the case when both $_A\bfX$ and
$\bfX_B$ are almost projective with   
$\Hom_{D^b(\Amod)}(\bfX,\bfX[1])=0$ and
$\Hom_{D^b(\modB)}(\bfX,\bfX[1])=0$.  If, in addition, $\chr K=0$, the first
Hochschild cohomology groups of $A$ and
$B$ are isomorphic as well.
\endproclaim

\demo{Proof} The first assertion is simply a translation  of Corollary 2.6 in
\cite{\GuSa} into the present context; the proof given in \cite{\GuSa} carries
over. The second part is an immediate consequence of Corollary 9 of the previous
section.  Since the Lie algebra associated to $\OA$ coincides with the first
Hochschild cohomology group of A in characteristic $0$, the last assertion
follows as well.
\qed\enddemo

\definition{Example} The above propositions cover situations which may be far
removed from the situation where $A$ and $B$ are derived equivalent. For
instance, if
$M$ is any nonzero left $A$-module with the properties that
$\Hom_A(M,A)=0$ and the $\Out$-orbit of $M$ is finite, then, by taking $X=A\oplus
M$ and
$B=\End(_AX)^{\op}\cong 
 \left[ \smallmatrix A & M\\ 0 & \End_A(M)^{\op} \endsmallmatrix \right]$
 we obtain a faithfully balanced bimodule $_AX_B$ satisfying the hypotheses of the
preceding propositions. Consequently, the centers of $A$ and $B$ are isomorphic,
as are
$\OA$ and $\OB$; the same is true for the first Hochschild cohomology groups of
$A$ and
$B$, provided that
$\chr K=0$.  But $A$ and $B$ fail to be derived equivalent, since they do not
have the same number of simple modules.  More specific examples of such modules
$M$ are the preprojective modules (in the sense of \cite{\AuSm}) which are devoid
of projective direct summands.  In case $\Out(A)$ is finite, {\it all} modules $M$
with vanishing $A$-dual $\Hom_A(M,A)$ qualify of course.  Finiteness of
$\Out(A)$, in turn, is guaranteed  whenever $A$ is a split algebra which is
tilting-cotilting equivalent to a hereditary algebra of tree type (cf.
\cite{\GuSa}, Theorem 2.10) \qed
\enddefinition

Our principal aim is to obtain the conclusion of Proposition 11 in case $A$ and
$B$ are derived equivalent algebras. For that purpose, we stretch the concept of a
faithfully balanced complex of bimodules as follows.

\definition{Definition} A  complex $_A\bfX_B$ of $A$-$B$-bimodules will be called
{\it derived faithfully balanced\/}, provided that the canonical $K$-algebra
homomorphisms $\lambda: A\longrightarrow \End_{D(\ModB)}(\bfX_B)$ and $\rho :
B\longrightarrow \End_{D(\AMod)}(_A\bfX)^{\op}$ are isomorphisms.  
\enddefinition 

\definition{Remark 12} Suppose that $_A\bfX_B$ is a derived faithfully balanced
bimodule complex.  Clearly, the canonical $K$-algebra homomorphism $\rho:
B\rightarrow
\End_{D(\AMod)}(_A\bfX)^{\op}$ equals the composition 
$$\rho: B\rightarrow \End_{C(\AMod)}(_A\bfX)^{\op}\rightarrow
\End_{D(\AMod)}(_A\bfX)^{\op},$$
 where the second map is the canonical one. If we set
$\Bhat =\End_{C(\AMod)}(_A\bfX)^{\op}$, our hypothesis on $_A \bfX_B$ thus yields
a vectorspace decomposition $\Bhat=B\oplus H$, where
$H$ is the two-sided ideal of $\Bhat$ consisting of the chain maps
$h:{}_A\bfX\rightarrow {}_A\bfX$ which become zero in
$D(\AMod)$. In the special case where $_A\bfX$ is an almost projective complex of
left
$A$-modules, these chain maps $h$ are precisely the ones which are homotopic to
zero.
\enddefinition

Of course, the preceding remark has a twin sibling, with the roles of $A$ and $B$
switched.  The same is true for the following useful observation.

\proclaim{Lemma 13} Let $_A\bfX_B$ be a bounded almost projective complex of
$A$-$B$-bimodules which is derived faithfully balanced. Then there is an
idempotent $e\in
\End_{\Cb}(_A\bfX)^{\op}$ with the following properties:

{\rm (1)} $_A\bfX(1-e)$ is an acyclic complex of left $A$-modules; in particular, 
$$B
\cong \End_{\Db}(_A\bfX e)^{\op}.$$

{\rm (2)} The kernel of the canonical $K$-algebra homomorphism
$$\End_{\Cb}(_A\bfX e) \rightarrow
\End_{\Db}(_A\bfX e)$$  is contained in the Jacobson radical of
$\End_{\Cb}(_A\bfX e)$.

{\rm (3)} A chain map in $\End_{\Cb}(_A\bfX e)$ is an isomorphism in $\Cb$ if and
only if it turns into an isomorphism in
$\Db$. 
\endproclaim

\demo{Proof} As in Remark 12, we denote by $\Bhat$ the endomorphism ring
$\End_{\Cb}(_A\bfX)^{\op}$, and by
$H$ the twosided ideal of those chain maps which become trivial in $\Db$.  Then
$\Bhat$ is in turn a finite dimensional algebra, because the complex $\bfX$ is
bounded by hypothesis.  The twosided ideal 
$\bigl(H+J(\Bhat)\bigr)/J(\Bhat)$ of $\Bhat/J(\Bhat)$ therefore gives rise to an
idempotent
$f\in H + J(\Bhat)$ such that $H+J(\Bhat)=f\Bhat+J(\Bhat)=\Bhat f+J(\Bhat)$.  We
verify that $f$ is trivial on the homology groups of $\bfX$. Indeed, if we write
$f$ in the form $f=h+j$, with $h\in H$ and
$j\in J(\Bhat)$, then $H^n(f) = H^n(j)$ for all $n$ by Remark 12. But $j$ being
nilpotent, so is the induced map $H^n(j)$, and hence $H^n(f) = 0$ due to
idempotency. Thus the chain complex $\bfX f$ is acyclic, and the choice $e=1-f$
satisfies condition (1).  Now the kernel $eHe$ of the canonical $K$-algebra
homomorphism $e\Bhat e\cong
\End_{\Cb}(_A\bfX e)^{\op}\rightarrow
\End_{\Db}(_A\bfX e)^{\op}\cong B$ is contained in $eJ(\Bhat)e$; for if it were
not, $H$ would fail to be contained in $f\Bhat+J(\Bhat)$, contrary to our choice
of $f$. To check the final point of our claim, we consider the canonical
surjection
$\End_{\Cb}(_A\bfX e)^{\op}\cong e\Bhat e\rightarrow e\Bhat e/eHe\rightarrow
e\Bhat e/eJ(\Bhat)e$. If $\beta\in
\End_{\Cb}(_A\bfX e)$ is a quasi-isomorphism, i.e., if $\beta$ becomes a unit in  
$e\Bhat e/eHe$, then, a fortiori, $\beta +eJ(\Bhat)e$ is a unit in
$e\Bhat e/eJ(\Bhat)e$.  But  
$J(e\Bhat e)=eJ(\Bhat)e$ being nilpotent, this implies that $\beta$ is a unit in
$e\Bhat e=\End_{\Cb}(_A\bfX e)^{\op}$, which proves (3).
\qed\enddemo

Our strategy in showing that the isomorphism type of $\OA$ is invariant under
derived equivalence is to focus on a certain closed subgroup of finite index of
$\Out(A)$; invariance of the identity component of this subgroup will then of
course amount to invariance of the identity component of the full outer
automorphism group.  More precisely, given mutually inverse twosided tilting
complexes
$_A\bfX_B$ and
$_B\Xtilde_A$ in the sense of
\cite{\Ric}, the subgroups of choice in $\Out(A)$ and $\Out(B)$ will be the
stabilizer subgroups
$\Out^{\bfX}(A)$ and $\Out^{\xtilde}(B)$ of $\bfX$ and $\Xtilde$ in $\Db$ and
$D^b(\Bmod)$, respectively.  The construction of an isomorphism between them
hinges on the following concept of a `semilinear chain map' between complexes of
left $A$-modules:  Fix
$\sigma \in \Aut(A)$ and complexes $\bfY, \bfZ$ of left $A$-modules.  Following
the pattern established in
\cite{\GuSa}, we call a chain map  $\varphi = (\varphi_n)_{n \in \NN}: \bfY
\rightarrow
\bfZ$ of $K$-complexes
$\sigma$-{\it semilinear}, provided that $\varphi_n(ay) = a^{\sigma}
\varphi_n(y)$ for all
$n \in \NN$, $y \in Y_n$, and $a \in A$.  We note that a $\sigma$-semilinear map
$\bfY
\rightarrow \bfZ$ can alternately be viewed as a morphism $^{\sigma^{-1}} \bfY
\rightarrow
\bfZ$ in $C(\AMod)$.  Moreover, we observe that, given any $\sigma$-semilinear
automorphism $\varphi$ of $\bfY$ and any morphism $\beta \in
\End_{C(\AMod)}(\bfY)^{\op}$, the conjugate $\varphi \beta \varphi^{-1}$ is again
a morphism in
$\End_{C(\AMod)}(\bfY)^{\op}$. To the endomorphisms (resp., automorphisms) of the
$K$-complex
$\bfY$ which are
$\sigma$-semilinear for some $\sigma \in \Aut(A)$ we collectively refer as the
($A$-){\it semilinear endomorphisms} (resp., {\it automorphisms}) of $\bfY$.  It
is straightforward that the set
$\S(\bfY)$ of all semilinear automorphisms of $\bfY$ forms a subgroup of the
group of automorphisms of the underlying $K$-complex.

In the following $\ltensor$ will denote the left derived functor of the tensor
functor on chain complexes.   Keep in mind that, under the hypothesis of the
following lemma, the algebra $e\Bhat e / eHe$ is canonically isomorphic to
$B$.        

\proclaim{Lemma 14} Again, let $_A\bfX_B$ be a bounded almost projective complex
of $A$-$B$-bimodules, which is derived faithfully balanced.  Moreover, let $e \in
\End_{\Cb}(_A\bfX)^{\op}$ be chosen as in Lemma 13, and let $H$ be the twosided
ideal of $\Bhat =
\End_{\Cb}(_A\bfX)^{\op}$ as specified in Remark 12. Given automorphisms
$\sigma\in
\Aut(A)$ and $\tau\in \Aut(B)$, conditions {\rm (1)}$-${\rm (3)} below are
related as follows: `$\roman{(1)}
\implies \roman{(2)}$' and `$\roman{(2)} \iff \roman{(3)}$'.

{\rm (1)} The functors $^\sigma (\bfX\ltensor_B-)$,
$\bfX\ltensor_B {}^{\tau}(-):D^b(\Bmod) \rightarrow D^b(\Amod)$ are isomorphic.

{\rm (2)} There is a $\sigma$-semilinear chain map $\xi :\bfX\rightarrow\bfX$,
i.e., a map in $\Hom_{\Cb}({} {}^{\sigma^{-1}}\bfX, \bfX)$, which turns into an
isomorphism
${}^{\sigma^{-1}}\bfX\rightarrow\bfX$ in $\Db$ and which satisfies the condition
that
$\xi b = b^{\tau} \xi$ in
$\Db$, for all $b\in B$.

{\rm (3)} There is a bijective $\sigma$-semilinear chain map $\varphi :\bfX
e\rightarrow\bfX e$ such that the automorphism of
$e\Bhat e$, given by
$\beta\mapsto \varphi\beta\varphi^{-1}$, induces the same map as $\tau$ on
$B$.
\endproclaim

\demo{Proof}
$(1)\implies (2)$.  One readily verifies that (1) is equivalent to the functors
${}^{\sigma^{-1}}(\bfX\ltensor_B -)$ and
$\bfX\ltensor_B{} ^{\tau^{-1}}(-)$ being isomorphic. We assume that this latter
condition is satisfied, and recall that, when applied to stalk complexes of
projectives, the functor $\bfX\ltensor_B -$ is just the usual tensor product (cf
\cite{\Ver}, ch\. 2). Now the tensor products  
$\bfX\otimes_B(^{\tau^{-1}}B)$ and $\bfX^\tau\otimes_BB$ are clearly canonically
isomorphic; we will identify them in fact.  Specifying an isomorphism   
$\chi :{}^{\sigma^{-1}}(\bfX\ltensor_B-) \rightarrow
\bfX\ltensor_B{}^{\tau^{-1}}(-)$ of functors $D^b(\Bmod) \rightarrow D^b(\Amod)$,
we consider $\chi(B) :{}^{\sigma^{-1}}(\bfX\otimes_B B)\rightarrow
\bfX\otimes_B{} ^{\tau^{-1}}B$ and the resulting chain of isomorphisms in
$D^b(\Amod)$: 
$${}^{\sigma^{-1}}\bfX\ \cong\ {}^{\sigma^{-1}}(\bfX\otimes_B B)\ \cong\ 
\bfX\otimes_B {}^{\tau^{-1}}B\ =\ \bfX^\tau\otimes_B B\cong \bfX^\tau.$$

Since ${}^{\sigma^{-1}}\bfX$ and $\bfX^\tau$ are both almost projective complexes
of left $A$-modules, the displayed isomorphism on the level of the derived
categories is induced by a chain map $\xi: {}^{\sigma^{-1}}\bfX\rightarrow
\bfX^\tau$ in $\Cb$ (see the remarks following the definition of an almost
projective complex in Section 1).  Alternately expressed, $\xi$ is a
$\sigma$-semilinear map
$\bfX \rightarrow \bfX^{\tau}$ which becomes bijective in $\Db$.  Finally, the
automorphism $\tau$ of $A$ makes a relevant appearance:  Namely, the naturality of
$\chi$ translates into the required additional property of $\xi$.

$(2)\implies (3)$. Let $\xi :\bfX\rightarrow\bfX$ be a
$\sigma$-semilinear map as specified in (2).  The fact that the morphism
$\bfX @>{e}>> \bfX e$ turns into an isomorphism in $\Db$ implies that, in the
derived category, $\xi$ can be viewed as an isomorphism
${}^{\sigma^{-1}}\bfX e\rightarrow\bfX e$, induced by the chain map
$\varphi = e\xi (e|_{\bfX e}) :{}^{\sigma^{-1}}\bfX e \rightarrow\bfX e$ in
$C^b(\Amod)$. From Lemma 13, part (3), one easily derives that
$\varphi$ is actually an isomorphism in $C^b(\Amod)$.  Thus $\varphi$ is a 
$\sigma$-semilinear automorphism $\bfX e\rightarrow\bfX e$. Checking the
remaining condition under (3) is routine. 

$(3)\implies (2)$. Using the equality $\bfX = \bfX e\oplus\bfX(1-e)$, one extends
$\varphi :\bfX e\rightarrow\bfX e$ to a $\sigma$-semilinear map
$\xi =\varphi\oplus 0:\bfX\rightarrow \bfX$ and notes that, trivially, the
extension has the required properties. \qed\enddemo

The following proposition represents the second crucial step  --  next to
Corollary 9  -- on the road towards the main result of this section.  We continue
to study a complex
$\bfX$ of $A$-$B$-bimodules as in Lemma 14, still denoting its endomorphism ring
in $\Cb$ by $\Bhat$.  Moreover, we introduce a set
$\G$ which will turn out to be the graph of a group isomorphism
$\phi: \Out^{\bfX}(A) \rightarrow \Out^{\xtilde}(B)$, provided that $_A\bfX_B$ and
$_B\Xtilde_A$ are mutually inverse twosided tilting complexes.  Namely, we define
$\G$ as the set of all pairs
$(\sigbar,\taubar)\in \Out(A)\times \Out(B)$ such that $\sigma$ and
$\tau$ satisfy the equivalent conditions (2), (3) of Lemma 14.  Since, from the
resulting definition of $\phi$, it is not obvious that this isomorphism of
abstract groups is actually an isomorphism of algebraic groups, we resort to
viewing $\phi$ from an alternate angle.  Namely, we will rewrite
$\phi$ in the form $p_2 p_1^{-1}$, with isomorphisms $p_1: \S \rightarrow
\Out^{\bfX}(A)$ and
$p_2: \S \rightarrow \Out^{\xtilde}(B)$ for a suitable choice of
$\S$, these latter maps having the benefit of being more readily recognizable as
isomorphisms of algebraic groups.  To that end, we consider the group
$\S(\bfX e)$ of semilinear automorphisms of the complex $\bfX e \in \Cb$,
together with the following two subgroups
$\U$ and $\V$:  The subgroup $\U$ consists of the maps of the form
$ x\mapsto uxv$, where $u$ is a unit of $A$ and $v$ a unit of $e \Bhat e$ (note
that the latter assignment is semilinear with respect to the automorphism $a
\mapsto uau^{-1}$ of $A$).  The subgroup 
$\V$ of $\S(\bfX e)$ consists of those maps
$\varphi\in\S(\bfX e)$ for which the algebra automorphism
$e\Bhat e\rightarrow e\Bhat e$, given by
$\beta\rightarrow\varphi\beta\varphi^{-1}$, induces an inner automorphism on
$e\Bhat e/eHe\cong B$; here again, $H \subseteq \Bhat$ is the ideal of those
chain maps which vanish in the derived category
$\Db$.   One readily checks that $\U$ and $\V$ are normal subgroups of $\S(\bfX
e)$ and that
$\U$ is contained in $\V$.  This provides us with a canonical projection 
$\pi:\S(\bfX e)/\U\rightarrow\S(\bfX e)/\V$.

\proclaim{Proposition 15} For any bounded almost projective complex $_A\bfX_B$ of
$A$-$B$-bimodules,  which is derived faithfully balanced, the following are true:

{\rm (1)} The set $\G$ is the graph of a group homomorphism $\phi
:\Out^{\bfX}(A)\rightarrow \Out(B)$, where $\Out^{\bfX}(A)$ is the
$\Out(A)$-stabilizer of $\bfX$ in $\Db$.

{\rm (2)} In case $\bfX$ is a tilting complex over $A$ and $e \in
\End_{\Cb}(\bfX)^{\op}$ an idempotent as before, the subgroups $\U$ and $\V$ of
$\S(\bfX e)$ are closed, as is the subgroup $\im(\phi)$ of
$\Out(B)$.  Moreover, there are isomorphisms of algebraic groups,
$$p_1:\S(\bfX e)/\U\rightarrow \Out^{\bfX}(A) \qquad\text{and}\qquad p_2:\S(\bfX
e)/\V\rightarrow \im(\phi ),$$  such that $p_2\pi =\phi p_1$, where $\pi$ is the
projection introduced above.
\endproclaim

\demo{Proof} (1) The task of checking that an element $\sigbar\in \Out(A)$ cannot
be the first component of two different pairs in $\G$ is easily reduced to the
case where $\sigbar=1$, i.e., to the case where $\sigma(a) = uau^{-1}$ for some
unit $u$ of $A$.  Suppose that $\tau \in \Aut(B)$ is such that the pair
$(\sigma,\tau)$ satisfies the equivalent conditions of Lemma 14, and let
$\varphi :\bfX e\rightarrow\bfX e$ be a bijective $\sigma$-semilinear
endomorphism of
$\bfX e$ as specified in condition (3) of that lemma; this is to say that
$\beta
\mapsto \varphi \beta \varphi^{-1}$ induces the same map as $\tau$ on $e \Bhat e/
e H e$.  The $\sigma$-semilinearity of $\varphi$ yields $A$-linearity of the
assignment
$v :\bfX e \rightarrow \bfX e$ defined by $xv = u^{-1}(x \varphi)$.  In other
words, $v$ is a unit in $\End_{\Cb}(\bfX e)^{\op} = e \Bhat e$. Keeping in mind
that the algebra
$e\Bhat e/ e H e$ is isomorphic to
$B$, by our balancedness hypothesis and the choice of $H$, we let $c$ be the
canonical image of $v$ in $B$.  In view of the equality $v\beta v^{-1} = \varphi
\beta
\varphi^{-1}$ for $\beta \in e \Bhat e$, we finally observe that $\tau$ is just
conjugation by
$c$, for our construction entails $\tau(b) = \varphi b \varphi^{-1} = c b c^{-1}$
for $b \in B$.  Thus $\tau$ is in turn inner, as required.      

Knowing that $\G$ is the graph of a function $\phi$, we pin down the domain of
$\phi$.  By construction, it consists of those classes $\sigbar\in
\Out(A)$ for which there exists a bijective $\sigma$-semilinear map
$\varphi: \bfX e\rightarrow\bfX e$ (note that conjugation by $\varphi$
automatically induces an automorphism of $e\Bhat e/ eHe$).  The latter condition
is tantamount to the requirement that
$\bfX e$ be isomorphic to ${}^\sigma (\bfX e)$ in $\Cb$, which is in turn
equivalent to the existence of an isomorphism $\bfX e\cong{}^\sigma (\bfX e)$ in
$\Db$;  this last equivalence is readily deduced from Lemma 13(3)  So the domain
of $\phi$ is  
$\Out^{\bfX}(A)$ as claimed. It is straightforward to check that $\phi$ is a group
homomorphism.

(2) Now suppose that $\bfX$ is a tilting complex of left A-modules.  In
particular, this means that the category $\add(\bfX)$, consisting of the finite
direct sums of direct summands of copies of $\bfX$, generates the homotopy
category $K^b(\Aproj)$ of all bounded projective complexes of finitely generated
left
$A$-modules as a triangulated category.  We infer that the left annihilator of
$\bfX$ in $A$ is zero.  Hence so is the annihilator of $\bfX e$, for the
projective complexes $\bfX$ and $\bfX e$ coincide in 
$K^b(\Aproj)$.  Following the model of
\cite{\GuSa, Lemma 1.2}, we deduce that, for every semilinear bijective chain map
$\varphi :\bfX e\rightarrow\bfX e$, there is a {\it unique\/} automorphism
$\sigma\in \Aut(A)$ with the property that $\varphi$ is $\sigma$-semilinear: 
Namely, if $\l_a$ denotes left multiplication of $\bfX e$ by $a$, then
$\sigma(a)$ is determined by the requirement that $\l_{\sigma(a)} = \varphi\l_a
\varphi^{-1}$.  If we write $\S$ for the group of semilinear automorphisms of
$\bfX e$, this provides us with a map
$\S\rightarrow
\Aut(A)$, which is actually a homomorphism of abstract groups. The proof of
assertion (1) shows that, whenever this homomorphism takes
$\varphi$ to an inner automorphism `$a
\mapsto uau^{-1}$' for a suitable unit $u \in A$,  the map
$v: \bfX e \rightarrow \bfX e$ sending $x$ to $u^{-1}(x\varphi)$ is a unit in
$e\Bhat e$, and vice versa.  But the latter condition just says that $x \varphi =
uxv$, and thus shows that our group homomorphism $\S \rightarrow \Aut(A)$ induces
a monomorphism  $\S/\U\rightarrow
\Out(A)$ of abstract groups.  In light of the previous paragraph, the image of
this monomorphism is $\Out^{\bfX}(A)$, and hence it gives rise to an isomorphism
$p_1:
\S/\U \rightarrow \Out^{\bfX}(A)$.    

Set $p_2' = \phi p_1 : \S/\U \rightarrow \Out(B)$.  Then, clearly $\im(p_2') =
\im(\phi)$.  As for the kernel of $p_2'$, we have $\phi p_1(\varphi +\U)=1$ if
and only if there is an inner automorphism
$\iota$ of $B$ such that $ b^\iota\varphi$ equals $\varphi b$ in
$\Db$, for all $b\in B$.  But this just means that the automorphism of $e\Bhat e$
given by $\beta\rightarrow
\varphi\beta\varphi^{-1}$ induces an inner automorphism of the algebra
$e\Bhat e/eHe\cong B$. Hence $\ker(p_2')$ equals $\V/\U$, so that $p_2'$ induces
an isomorphism $\S/\V \rightarrow \im(\phi)$, as required.

Finally, we need to ascertain that $p_1, p_2$ are isomorphisms of algebraic
groups.  As for $p_1$, we start by identifying the algebra $A$ with the
subalgebra of $\Ahat = \End_{C^b(\moddash e {\hat B} e)}(\bfX e)$ consisting of
the left multiplications $\l_a$ on $\bfX e$ for $a\in A$; this identification is
legitimate, since
$\bfX e$ is faithful over $A$.  With $A$ viewed from this angle, the group
homomorphism
$\S \rightarrow \Aut(A)$ underlying
$p_1$ maps $\varphi \in \S$ to the automorphism $\l_a \mapsto \varphi \l_a
\varphi^{-1}$ of $A$, the latter conjugate being equal to $\l_{\sigma(a)}$ if
$\varphi$ is $\sigma$-semilinear.  This assignment is clearly a morphism of
algebraic groups, and hence $\U$ is a closed (normal) subgroup of $\S$.  To see
that the inverse of $p_1$ is again a morphism, it suffices to show that the map
$p_1': \S \rightarrow \Out^{\bfX}(A)$ which induces $p_1$ modulo $\U$ has a
surjective differential $\T_e(\S) \rightarrow
\T_e(\im(p_1'))$ (cf. \cite{\Bor, Proposition 6.13}).  

In a first step, we focus on the following morphism of algebraic groups
$$\gamma: G \rightarrow \GL(V) \qquad \text{with} \qquad \gamma_g(v) = gvg^{-1},$$
where $G = \GL(\bfX e)$ and $V = \End_K(\bfX e)$. Our aim is to see that the
restriction $\S
\rightarrow \gamma(\S)$, again denoted by $\gamma$, has surjective differential
$\d\gamma: \T_e(\S) \rightarrow \T_e(\im(\gamma)) \subseteq \End_K(V)$.  From
\cite{\Bor, 3.10}, we know that
$$ \qquad (\d\gamma)(X)(v) = (\d \frak{orb}_v)_e(X), \tag{$\dagger$}$$  for $X
\in \T_e(\S)$ and $v \in V$, where $\frak{orb}_v: \S \rightarrow \gamma(\S)v$ is
the orbit map $g \mapsto \gamma_g v$.   To check surjectivity, we let $v_1,
\dots, v_n$ be a basis for
$V$, set
$\underline{v} = (v_1, \dots, v_n) \in V^n$,  and consider the `expanded orbit
map'
$$\Orb: G \rightarrow \gamma(G).\underline{v}, \qquad g \mapsto g.\underline{v} =
(\gamma_g v_1, \dots, \gamma_g v_n),$$ and its restriction to $\S$,
$$\orb: \S \rightarrow \gamma(\S).\underline{v}.$$ The map $\Orb$ has reduced
fibres, since the stabilizer subgroup
$\Stab(\underline{v})$ of
$\underline{v}$ in $G$ clearly arises as the solution set of a system of linear
equations (for background on reduced fibres, see \cite{\Kra, A.I.2.5, 2.6}).  By
\cite{\Kra, A.I.5.5, Satz 2}, this implies that $\ker(\d \Orb)_e =
\T_{\underline{v}}(\Stab(\underline{v}))$.  Cutting down $\Stab(\underline{v})$ to
the stabilizer subgroup $\stab(\underline{v})$ of $\underline{v}$ in $\S$, we
infer that
$\ker(\d \orb)_e$ equals $\T_{\underline{v}}(\stab(\underline{v}))$, and, in view
of
\cite{\Bor, Proposition 6.7}, we deduce further that the differential $(\d
\orb)_e:
\T_e(\S)
\rightarrow \T_{\underline{v}}(\gamma(\S).\underline{v})$ is surjective.  In
light of
$(\dagger)$, this yields surjectivity of $\d\gamma$.  

If we enlarge the codomain of $p_1'$ from $\Out^{\bfX}(A)$ to $\Out(A)$, then
$p_1'$ can be factored in the form
$$\S @>\gamma>> \gamma(\S) @>\res>> \Aut(A) @>\can>> \Out(A),$$ where $\res:
\gamma(\S) \rightarrow \Aut(A)$ sends any $K$-automorphism of $V$ in
$\gamma(\S)$ to its restriction to the invariant subspace $A = \Ahat$.  Now
$\res: 
\gamma(\S) \rightarrow \im(\res)$ is clearly separable, and so is the canonical
map $\can$, the latter being actually a geometric quotient, whence both  $\res$
and
$\can$ have surjective differentials.  This shows that our auxiliary map $p_1'$
has surjective differential, and we conclude that $p_1$ is indeed an isomorphism
of algebraic groups.     

To deal with $p_2$, we identify $B$ with $e\Bhat e/eHe$ via the canonical
isomorphism.  Clearly, the group homomorphism $\S \rightarrow \Aut(e
\Bhat e)$ sending
$\varphi$ to the automorphism $\beta \mapsto \varphi \beta \varphi^{-1}$ respects
the variety structures, and since the normal subgroup of $\Aut(e \Bhat e)$
consisting of those automorphisms which induce inner automorphisms of $e \Bhat
e/eHe$ is closed, this subgroup  gives rise to a geometric quotient of $\Aut(e
\Bhat e)$.  In particular, the auxiliary map $p_2':\S / \U \rightarrow \Out(B)$
as above is again a morphism of algebraic groups.  This in turn entails that  $\V
\subseteq \S$ and $\im(\phi) \subseteq \Out(B)$ are closed subgroups. Now the
geometric quotient property of the canonical map $\S \rightarrow
\S/\V$ guarantees that $p_2$ is a morphism of algebraic groups as well.  With an
argument similar to that given for $p_1$, one finally shows that $p_2$ also has a
surjective differential, which guarantees that its inverse is in turn a morphism
of algebraic groups.
\qed\enddemo

In the following, we  adopt Rickard's terminology in \cite{\Ric}, but
insignificantly deviate from his conventions by calling a category equivalence
$D^-(\BMod)
\rightarrow D^-(\AMod)$ (resp., $D^b(\Bmod) \rightarrow D^b(\Amod)$) {\it
standard\/} if it is isomorphic to $\bfX \ltensor_B - $, where $_A \bfX_B$ is a
two-sided tilting complex; as is shown in \cite{\Ric, Theorem 3.3}, this is
harmless. Moreover, we note that   according to
\cite{\Ric, loc. cit.}, any equivalence
$\bfX
\ltensor_B -$ between the derived categories $D^-(\BMod)$ and $D^-(\AMod)$
automatically restricts to an equivalence $D^b(\Bmod) \rightarrow D^b(\Amod)$. 
To finally clear the road to the main theorem of this section, we record an easy
subsidiary observation.  

\proclaim{Lemma 16} Let B be a finite dimensional algebra and
$F:D^-(\BMod)\rightarrow D^-(\BMod)$ a standard derived equivalence. Then the
following assertions are equivalent:

{\rm (1)} $F$ is induced by an element $\taubar\in \Out(B)$, i.e., there is a
natural isomorphism $F\cong {}^\tau (-)$.

{\rm (2)} $F$ takes the stalk complex $_B B$  to an isomorphic copy of itself.
\endproclaim

\demo{Proof} Clearly, any $\tau\in \Aut(B)$ gives rise to an isomorphism $B\cong
{}^{\tau}B$ of left
$B$-modules, which proves `$(1) \implies (2)$'.  For the converse, let 
$_B\bfT_B$ be a two-sided tilting complex for $B$ such that $F\cong
\bfT\ltensor_B -$. Then
$B\cong F(B)\cong \bfT\ltensor_B B\cong \bfT$ in $D^-(\BMod)$.  In other words,
as a complex of left $B$-modules, $\bfT$ is isomorphic to the stalk complex $B$
on the level of the derived category.  By \cite{\Ric}, the right derived functor
of $\Hom_B(\bfT,-)$ is a quasi-inverse of $F$ -- call it  $F^{-1}$. We deduce
that $F^{-1} =  R\Hom_B(\bfT,-)\cong
\Hom_{D^-(\BMod)}(B,-)$ is induced by a Morita equivalence
$\Hom_B(C, -): \BMod\rightarrow \BMod$, where $C$ is some $B$-$B$-bimodule
isomorphic to
${}^1B^{\tau^{-1}}$ for some automorphism
$\tau$ of $B$.  This shows that $F \cong {}^{\tau}(-)$ is determined by $\taubar
\in
\Out(B)$, and completes the proof of `$(2) \implies (1)$'. \qed\enddemo

As announced, we will exploit the scenario we have rigged up in the special case
where $_A
\bfX_B$ is a twosided tilting complex.  In particular, we will see that, in this
situation, the map
$\phi$ constructed in Proposition 15 is a group isomorphism $\Out^{\bfX}(A)
\rightarrow
\Out^{\xtilde}(B)$, provided that $_B\Xtilde_A$ is a twosided tilting complex
inverse to
$\bfX$.  Moreover, we will obtain $\U = \V$, which will provide us with an
alternate description of $\phi$ as $p_2 p_1^{-1}$, in the terminology of the
proposition.

\proclaim{Theorem 17} Suppose that $A$ and $B$ are finite dimensional algebras
(over an algebraically closed field
$K$).  If
$A$ and $B$ are derived equivalent, then $\OA$ and $\OB$ are isomorphic algebraic
groups. 
\endproclaim 

\demo{Proof} 
  We begin by choosing two-sided tilting complexes
$_A\bfX_B$ and
$_B\Xtilde_A$ which are inverse to each other, that is, $\bfX \ltensor_B
\Xtilde \cong {} _A A_A$ and
$\Xtilde \ltensor_A \bfX \cong {}_B B_B$ in the derived categories $D^b(A\otimes_K
A^{\op}\dashmod)$ and $D^b(B\otimes_K B^{\op}\dashmod)$, respectively.  Recall
that
$\Out^{\bfX}(A)$ and
$\Out^{\xtilde}(B)$ denote the subgroups consisting of those elements of $\Out(A)$
and
$\Out(B)$ which stabilize $\bfX$ and
$\Xtilde$ in $\Db$ and $D^b(\Bmod)$, respectively.  Our first goal is to show
that these closed subgroups of $\Out(A)$ and $\Out(B)$ have finite index in the
corresponding full outer automorphism groups.  

We show finiteness of the index $[\Out(A): \Out^{\bfX}(A)]$, the result for $B$
being symmetric.  Note that the $\Out(A)$-stabilizer of $\bfX$ in $\Db$ equals
that of $\bfX e$ in $\Db$, where $e$ is an idempotent in the $\Cb$-endomorphism
ring of $\bfX$ as specified in Lemma 13, for, in the derived category, $\bfX$ and
$\bfX e$ become isomorphic.  Hence the index under consideration coincides with
$[\Out(A):
\Out^{\bfX e}(A)]$ and thus equals the cardinality of the $\Out(A)$-orbit of
$\bfX e$ in $\Db$.  So clearly it suffices to show that the $\Out(A)$-orbit of
$\bfX e$ in $\Cb$ is finite.  But this follows from Corollary 9 of Section 1,
which shows that there are only finitely many complexes over $A$ which share the
$K$-dimension of
$\bfX e$, as well as the property that $\Hom_{\Db}(\bfX e, \bfX e[1]) = 0$.  This
takes us to our first goal.  

Our argument  is thus reduced to showing that $\im(\phi) = \Out^{\xtilde}(B)$
and  $\U = \V$, in the terminology of Proposition 15.  Once these equalities are
established, Proposition 15 will provide us with isomorphisms of algebraic groups
$$\Out^{\bfX}(A) @<p_1<< \S(\bfX e)/\U @>p_2>> \Out^{\xtilde}(B).$$ We deduce
that $\phi= p_2 p_1^{-1}$ restricts to an isomorphism $\OA
\rightarrow \OB$, for the identity component of $\S(\bfX e)/\U$ is mapped onto
that of
$\Out^{\bfX}(A)$  --  which equals $\OA$  --  by $p_1$, and onto that of
$\Out^{\xtilde}(B)$  --  which equals $\OB$  --  by
$p_2$.  

We access the two remaining equalities by way of the group homomorphism $\phi$
from
$\Out^{\bfX}(A)$ to $\Out(B)$ constructed in Proposition 15.  To that end, we
give yet another  description of $\phi$ in the present situation.  For any element
$\sigbar \in
\Out^{\bfX}(A)$, the equivalence
$F = (\Xtilde \ltensor_A -)\circ {}^{\sigma}(\bfX \ltensor_B -)$ from
$D^b(\Bmod)$ to
$D^b(\Bmod)$ preserves the isomorphism type of the stalk complex $_B B$, whence
Lemma 16 shows $F$ to be induced by some element $\taubar \in \Out(B)$.  In view
of the choice of
$\bfX$ and
$\Xtilde$ as inverse to each other, this provides us with an isomorphism
${}^\sigma (\bfX\ltensor_B  -)\cong\bfX\ltensor_B {}^\tau(-)$ of functors
$D^b(\Bmod) \rightarrow
\Db$.  Consequently, Lemma 14 guarantees that the pair
$(\sigbar, \taubar)$ satisfies the equivalent conditions (2), (3) of that lemma;
in other words, the pair $(\sigbar, \taubar)$ belongs to the graph $\G$ of
$\phi$.  Thus
$\phi$ assigns to any $\sigbar \in
\Out^{\bfX}(A)$ the unique element $\taubar \in \Out(B)$ with the property that
${}^\sigma (\bfX\ltensor_B  -)\cong\bfX\ltensor_B {}^\tau (-)$.  To obtain
$\im(\phi)
\subseteq
\Out^{\xtilde}(B)$, we  derive the following string of isomorphisms in
$D^b(\Bmod)$:
$${}^\tau\Xtilde\cong (\Xtilde \ltensor_A \bfX)\ltensor_B {}^\tau
\Xtilde \cong\Xtilde\ltensor_A  (\bfX\ltensor_B
{}^\tau\Xtilde)\cong\Xtilde\ltensor_A {}^\sigma (\bfX\ltensor_B
\Xtilde)\cong\Xtilde\ltensor_A {}^\sigma A.$$ Since the left $A$-modules
${}^\sigma A$ and $A$ are isomorphic,  we conclude that ${}^{\tau}\Xtilde \cong
\Xtilde$ in $D^b(\Bmod)$.  But this just means that $\taubar \in 
\Out^{\xtilde}(B)$.

The symmetric version of Proposition 15, obtained by switching the roles of $A$
and $B$, permits us to `symmetrize' the preceding paragraph.  Indeed, it provides
us with a group homomorphism
$\widetilde{\phi}:\Out^{\xtilde}(B)\rightarrow \Out(A)$ which assigns to $\taubar
\in
\Out^{\xtilde}(B)$ the unique element
$\sigbar\in \Out(A)$ with ${}^\tau(\Xtilde \ltensor_A  -)\cong\Xtilde\ltensor_A
{}^\sigma (-)$, and we obtain $\im(\widetilde{\phi})\subseteq
\Out^{\bfX}(A)$ as above. From the fact that
$\bfX$ and  $\Xtilde$ are inverse to  each other, we moreover derive that $\phi$
and $\widetilde{\phi}$ are inverse group isomorphisms between $\Out^{\bfX}(A)$
and $\Out^{\xtilde}(B)$. In particular, $\im(\phi) = \Out^{\xtilde}(B)$.  In view
of part (2) of Proposition 15, finally, injectivity of
$\phi$ guarantees that the canonical map $\pi: \S(\bfX e)/\U \rightarrow \S(\bfX
e)/\V$ is the identity; in other words, we also have $\U = \V$.  This completes
the argument. 
\qed \enddemo

\head Acknowledgement\endhead

The paper was prepared while the second-named author was visiting the University
of California at Santa Barbara. He would like to thank the UCSB Mathematics
Department for its hospitality. Moreover, he is grateful to the D.G.E.S. of Spain
and the Fundaci\'on `S\'eneca' of Murcia for their financial support.

\enddocument